\documentclass[10pt,twoside]{article}
\usepackage{amsmath, amsfonts, amsthm, amssymb, amscd}
\usepackage[mathcal]{euscript} 
\usepackage{dsfont, pxfonts, verbatim}
\usepackage{mathrsfs}  

\usepackage{a4wide}

\newtheorem{theorem}{Theorem}[section]
 
\newtheorem{lemma}[theorem]{Lemma}
\newtheorem{corollary}[theorem]{Corollary} 
\newtheorem{definition}[theorem]{Definition} 
  
\numberwithin{equation}{section} 
\newcommand \boundd {C_0} 
\newcommand \param {\delta}
\newcommand \alphab {\overline \alpha} 
\newcommand \be           {\begin{equation}}
\newcommand \ee            {\end{equation}}

\newcommand \RR           {\mathbb{R}}

\newcommand \gb             {{\overline g}}
\newcommand \gbar {{\overline g}}

\newcommand \del           \partial
\newcommand \eps            \epsilon
\newcommand \ubar       {{\overline u}}
\newcommand \cbar       {{\overline c}}
\newcommand \barc       {{\underline c}}

\newcommand \loc        {{\mathrm{loc}}}

\DeclareMathOperator    \sgn {sgn}
\DeclareMathOperator    \TV  {TV}

\DeclareMathOperator\Lip {Lip}

\newcommand \la         \langle
\newcommand \ra         \rangle

\newcommand{\clg}[1]{{\mathcal{#1}}}

\newcommand \Ocal   {\mathcal O} 
\newcommand{\Hcal}{\mathcal{H}} 
\newcommand{\Lorentz}{{M}} 
\newcommand{\Lorentzplus}{{M}_+}
 
\DeclareMathOperator\divex {div}
\newcommand{\dive}{\divex_{g}} 
\DeclareMathOperator \supp {supp} 

\newcommand{\bfOmega} \Omega 

\newcommand \alphabar {\overline \alpha}
\def\Xint#1{\mathchoice
{\XXint\displaystyle\textstyle{#1}}%
{\XXint\textstyle\scriptstyle{#1}}%
{\XXint\scriptstyle\scriptscriptstyle{#1}}%
{\XXint\scriptscriptstyle\scriptscriptstyle{#1}}%
\!\int}
\def\XXint#1#2#3{{\setbox0=\hbox{$#1{#2#3}{\int}$}
\vcenter{\hbox{$#2#3$}}\kern-.5\wd0}}

\def\dashint{\Xint\diagup}

\newcommand{\one}{\mathds{1}}

\begin{document}
\title{
A geometric approach to error estimates for conservation laws posed on a spacetime
}
\author{Paulo Amorim$^1$, Philippe G. L{\large e}Floch$^2$, and Wladimir Neves$^3$}

\date{}

\maketitle

\footnotetext[1]{Centro de Matem\'atica e Aplica\c c\~oes
Fundamentais, Departamento de Matem\'atica, Universidade de Lisboa, Av.~Prof.~Gama~Pinto 2,
1649-003 Lisboa, Portugal. 
E-mail: {\sl pAmorim@ptmat.fc.ul.pt}}

\footnotetext[2]{Laboratoire Jacques-Louis Lions \& Centre
National de la Recherche Scientifique, Universit\'e Pierre et Marie Curie (Paris 6), 
4 Place Jussieu,  75252 Paris, France. E-mail: {\sl pgLeFloch@gmail.com}}

\footnotetext[3]{Instituto de Matem\'atica, Universidade Federal
do Rio de Janeiro, C.P. 68530, Cidade Universit\'aria 21945-970,
Rio de Janeiro, Brazil. E-mail: {\sl Wladimir@im.ufrj.br}
\newline
%\noindent
 Revised in April 2011. To appear in: Nonlinear Analysis (2011). 
}

\begin{abstract} We consider a hyperbolic conservation law posed on an 
$(N+1)$-dimensio\-nal spacetime,
whose flux is a field of differential forms of degree $N$. 
Generalizing the classical Kuznetsov's method,
we derive an $L^1$ error estimate which applies
to a large class of approximate solutions. 
In particular, we apply our main theorem and deal with 
two entropy solutions associated with distinct flux fields, as well as 
with an entropy solution and an approximate solution.  Our framework 
encompasses, for instance, equations posed on a globally hyperbolic Lorentzian manifold. 
\end{abstract}

%\tableofcontents

%=====================================================================

\section{Introduction}
 
This paper provides a general
framework leading to error estimates for hyperbolic conservation laws
posed on an $(N+1)$-dimensional manifold $M$, referred to as a spacetime and, in particular, 
leading to a sharp estimate for the difference, measured in the $L^1$ norm, between an exact solution 
and an approximate solution. The present paper can be regarded as a generalization to manifolds
of a contribution by Bouchut and Perthame \cite{BP}, who recast in a concise form the
pioneering works of Kruzkov and Kuznetsov \cite{K,Kuz1,Kuz2} for hyperbolic conservation
laws posed on the flat (Euclidian) spacetime. We are thus interested here in extending these results 
to conservation laws defined on manifolds, and develop a physically more realistic setting when geometrical effects are now taken into account.

Motivated by the case of the shallow water equations on the sphere, 
the theory of hyperbolic conservation laws on manifolds has been developed in recent years by 
LeFloch together with collaborators. In particular, well-posedness results have been obtained in 
\cite{BL,LeFloch,LO}, and convergence results for finite volume schemes in \cite{ABL,ALO}, while an error 
estimate in the case of a Riemannian manifold was derived in \cite{LNO}. 
For further results on the well-posedness theory, we also refer the reader to contributions by 
Panov in~\cite{Panov1,Panov2} and, on the finite volume schemes, to the earlier work \cite{CCL}.

Recently, in \cite{LO}, LeFloch and Okutmustur introduced a framework based on differential forms 
and dealt with conservation laws defined on an $(N+1)$-dimensional manifold $M$. 
In their formulation, the flux of the equation is given by a field of $N$-forms, 
rather than by a vector field as was the case in earlier works. The formulation based on $N$-forms 
is geometrically natural in that only minimal assumptions on the geometrical structure are 
required in order to establish the well-posedness of the initial value problem.

Considering hyperbolic conservation laws posed on a spacetime,  we derive here a coordinate-free
error estimate in the $L^1$-norm, which involves, on one hand,
 an exact weak solution satisfying all entropy inequalities
and, on the other hand, an approximate solution for which the entropy dissipation measures contain
 some
(positive) ``error terms''.  Our assumptions on 
the manifold are minimal and, in particular, no metric is needed and no volume form is a priori prescribed. 
Indeed, our method only requires a suitable ``hyperbolicity condition'', 
and the prescription of a certain family of mollifiers (similar to mollifiers that may be determined from the distance function when a Riemann metric is given on the manifold). 
Still, we observe below that, to the flux-field defining the conservation law, one can 
associate a natural choice of volume form defined on the manifold. 
In addition, from our main theorem 
we deduce, on one hand, 
the $L^1$-contraction property of the semigroup of entropy solutions 
and, on the other hand,
a sharp estimate comparing together the solutions of two conservation laws with different flux fields.

In the second part of this paper, we investigate the case 
of a conservation law defined on a Lo\-rentz\-ian manifold. Here, the presence of 
a metric allows us to refine our estimates. For instance, a suitable family of mollifiers is provided by the 
metric structure of the manifold, and the error terms can be given in a more explicit way, since they are naturally 
written using the metric. More importantly, the introduction of a metric allows us to consider second-order
error terms in the approximate solutions. This allows us to apply our main 
result further, and  derive error estimates for 
a nonlinear diffusion model.

We may summarize the main difficulties overcome in this paper as follows. 
First of all, since the conservation law under consideration is posed on a non-flat
manifold, the geometry of that manifold must be taken
into account in, for instance, the formulation of approximate schemes and the analysis of their convergence. 
It has been pointed out that geometric effects occur 
which change the qualitative properties of solutions \cite{ABL}. 
Spacetimes and, in particular, Lorentzian manifolds may not be ``invariant by translation'' in the time direction, 
and the time and spatial geometries are intertwined, giving rise 
to phenomena not present in the Euclidian or Riemannian set-up.  
The main difficulty dealt with here 
lies in the lack of geometric structure on the spacetime which 
makes it difficult to 
deploy the analytical techniques used when a metric is prescribed. 
To circumvent this problem, we assume the existence of a suitable family of mollifiers, which 
are adapted to the sole structure available on the manifold, namely the family of flux fields of the 
conservation law. For this reason, our estimate in Theorem~\ref{DF-10} is, later on in 
 Theorem \ref{10-10}, specialized to the case that a metric is specified.
 (See also an earlier result in \cite{LNO}.) 
 Importantly, all of the estimates established in the present paper are 
coordinate-free.

An outline of the paper follows. Section~\ref{CL} is devoted to the general framework.
First of all, in Section~\ref{21s}, we recall 
from \cite{LO} some basic concepts about conservation laws posed on a spacetime. We then 
introduce the notion of admissible mollifiers, which allows us to state our main result, in 
Theorem \ref{DF-10} below. Then, in Section~\ref{24s}, 
we discuss our first application concerning two conservation laws with distinct flux fields.
The rest of the section is devoted to the proof of Theorem~\ref{DF-10}. 
In the second part of the paper encompassing the whole of Section~\ref{LM}, 
we treat the special case of conservation laws posed on a Lorentzian manifold, and 
state our error estimate in Theorem~\ref{10-10}. 
We are then in a position to provide two more applications in Section~\ref{33s}
and \ref{ND}. Section~\ref{230} contains a proof of Theorem~\ref{10-10}.

%=====================================================================

\section{Error estimates for a spacetime}
\label{CL}

\subsection{Conservation laws based on differential forms}
\label{21s}

Let $M$ be an $(N+1)$-dimensional manifold (with smooth topological structure), which we refer to as a \emph{spacetime}.
Denote by $\Lambda^k(M)$ the space of all smooth fields of $k$-differential forms on $M$, 
and by $d: \Lambda^k(M) \to \Lambda^{k+1}(M)$ the exterior derivative operator. 

Following LeFloch and Okutmustur \cite{LO}, we recall the formulation of the initial value problem 
for hyperbolic conservation laws posed on a spacetime. The setting is based on differential forms
and the conservation law reads 
\be 
\label{DF.10}
 d(\omega(u)) = 0 \qquad \text{on} \quad M,
\ee
where $u: M \to \RR$ is the unknown function and 
the given family of smooth differential $N$-forms $\omega=\omega(u) \in \Lambda^N(M)$ 
depends smoothly upon the variable $u$ and is referred to as the {\sl flux field} 
of the conservation law. 
The field $\omega$ is said to be {\sl geometry compatible} if 
it is exact, that is,
\[
 (d\omega)(\ubar) = 0, \qquad  \ubar \in \RR.
 \]
Moreover, \eqref{DF.10} is supplemented with the initial 
condition 
\be \label{DF.101}
   u|_{\clg{H}_0} = u_0,
\ee
where $u_0$ is a given data
defined on a hypersurface $\clg{H}_0 \subset M$. Throughout this paper we assume that the data 
$u_0$
and, therefore, the unknown function (thanks to the maximum principle)
is bounded, that is, 
\be
\label{boundM}
-\boundd \leq u_0 \leq  \boundd,
\ee
for some $\boundd>0$. We tacitly restrict all values of $u$ to 
lie in the interval $[-\boundd,\boundd]$, and 
we point out that 
all of our estimates in this paper depend implicitly on this constant $\boundd$.

Following \cite{LO}, we introduce the following notion of global hyperbolicity:
the manifold $M$ is foliated by hypersurfaces, 
\be \label{MHT}
M = \bigcup_{t\geq 0} \clg{H}_t 
\ee
for some (smooth) time-function $t$, where each slice $\Hcal_t$ has the same topology as $\Hcal_0$ which we assume to be a (smooth), compact, orientable $N$-manifold. 
For each $T>0$, we use the notation
$$
M_T := \bigcup_{0 \leq t \leq T} \clg{H}_t. 
$$
Moreover, given any hypersurface $\Hcal$, we denote by $i : \Hcal \to M$ the canonical imbedding 
and by $i^*$ its pullback, taking forms in $\Lambda^k(M)$ to forms in $\Lambda^k(\Hcal)$.

\begin{definition}
\label{NDG}  
Let $M$ be a foliated manifold \eqref{MHT} and $\omega=\omega(\ubar)$ be a flux field. 
The conservation law \eqref{DF.10} is said to satisfy the 
{\bf global hyperbolicity condition}
if for all $t \geq 0$ the $N$-form field $i^*_{\Hcal_t} \omega(0)$ is a volume form on $\Hcal_t$, and 
there exist constants $0< \barc <   \cbar$ independent of $t$ such that for all $\ubar \in [-C_0, C_0]$ 
\be 
\label{GHC}
\barc \, i^*_{\Hcal_t} \del_u \omega(0) \leq  \, i^*_{\Hcal_t} \del_u \omega(\ubar) \leq  \cbar  \, i^*_{\Hcal_t} \del_u \omega(0), 
\ee
as inequalities between $N$-forms defined on the slice $\Hcal_t$. 
\end{definition}

In particular, the condition \eqref{GHC} implies that for each non-empty smooth hypersurface $e \subset \Hcal_t$
the integral
$$
\int_e i_e^* \partial_u \omega(\ubar)
$$
is positive, and its ratio with $\int_e i_e^* \partial_u \omega(0)$ is bounded above and below. 

Recall also the following definition of entropy solutions.

\begin{definition}
A smooth field of $N$-forms $\bfOmega = \bfOmega(\ubar)$ is called a 
{\bf (convex) entropy flux field}
 for the conservation law \eqref{DF.10} if there exists a (convex) function $U :\RR \to \RR$ such that
\[
\bfOmega(\ubar) = \int_0^\ubar \del_u U(v) \del_u \omega(v) \, dv, \qquad \ubar \in \RR.
\]
A measurable and bounded function $u: M\to \RR$ is called an
{\bf entropy solution} to the 
Cauchy problem \eqref{DF.10}-\eqref{DF.101} if for every
(compactly supported and smooth) non-negative test-function $\varphi$,
the following {\bf entropy inequalities} hold for any convex {\bf entropy pair} $(U, \Omega)$ 
\be 
\label{DF.20} 
\aligned
\int_{M} \Big( d\varphi \wedge \bfOmega(u) + \varphi (d \bfOmega)(u) 
- \varphi \del_u U(u) (d\omega)(u) \Big) + \int_{\Hcal_0} \varphi|_{\Hcal_0} i_{\Hcal_0} ^*\bfOmega(u_0) \ge 0.
\endaligned
\ee
\end{definition}
 
Under the global hyperbolicity condition above, the flux of the equation naturally induces an  
($(N+1)$-dimensional) {\sl reference volume form} on $M$, namely
\be 
\label{DF.65} 
\alphabar := dt\wedge \del_u\omega (0). 
\ee
We emphasize that this volume form depends on, both, the flux field (at the state $\ubar=0$) and
the chosen foliation. 
On the other hand, a more fundamental structure on $M$ is provided by the family of {\sl $N$-form} flux  
 $\omega=\omega(\ubar)$
which determines the conservation law under consideration. 

Finally, we recall from \cite{LO} that, under the global hyperbolicity condition,
the initial value problem \eqref{DF.10}-\eqref{DF.101} admits a unique entropy solution
that depends Lipschitz continuously in the $L^1$ norm upon its initial data.

%----------------------------------------------------------------------------------------

\subsection{Approximate solutions}  

The main objective of this paper is to provide a general framework for the derivation of 
error estimates for hyperbolic conservation laws defined on a spacetime. 
In our statements and proofs we will use the \textsl{Kruzkov entropies}
\[
U(u,k) = | u - k|, \qquad k \in \RR
\]
with entropy flux
\[
\Omega (u,k) = \sgn(u - k) (\omega (u) - \omega(k)).
\]
Let $v$ be an (exact) entropy solution to the
conservation law \eqref{DF.20}, satisfying therefore 
\be 
\label{DF.60}
   d(\bfOmega(v,l) ) + G(v,l) \le 0,  \qquad l \in \RR 
\ee 
in the sense of distributions, 
where $G= G(v,l)$ are $(N+1)$-forms defined by
\[
G(v,l) := \sgn(v - l) (d\omega) (l).
\]

Following now Bouchut and Perthame~\cite{BP} and in order to deal with approximate solutions, it is convenient to introduce Radon measures
to estimate error terms that arise in an approximate version of the entropy inequalities. 
Furthermore, in our setting, rather than functions or vector fields, 
we have to deal with $N$-form fields on the manifold $M$. 
Hence, we now introduce Radon measure-valued fields of $N$-forms. 
For instance, a distributional $(N+1)$-form is an element of 
the dual of the space of test-functions on $M$, while a (scalar) distribution is an element of the dual of the space 
of (smooth) compactly supported $(N+1)$-forms. 
In what follows, to keep the notation simple we write, for instance, 
\[
\int_M \varphi \, \rho 
\]
for the duality bracket between a Radon-measure $(N+1)$-form field $\rho$ 
and a (continuous) test-function $\varphi$.

We are now in a position to write the {\sl approximate entropy inequalities} 
satisfied by some approximate solution $u: M \to \RR$, that is,  
\be 
\label{DF.40} \aligned d(\bfOmega(u,k) ) + G(u,k) \le d H_k  +
K_k =: E_k,  \qquad k \in \RR, 
\endaligned
\ee
where $H_k$ is a family of locally Radon measure-valued $N$-form fields and 
$K_k$ a family of locally Radon measure-valued $(N+1)$-form fields. We make the following key assumption
on these error terms: 
there exist a non-negative Radon measure-valued $N$-form field $\alpha_H$
and a non-negative Radon measure-valued  $(N+1)$-form field $\alpha_K$ such that 
for every $1$-form test-field $\gamma$ and every test-function $\varphi \geq 0$ 
\be 
\label{DF.50} 
\aligned
& \sup_{k \in \RR} \Big| \int_M H_k \wedge \gamma \Big| \leq \Big| \int_M \alpha_H \wedge\gamma \Big|, 
\\
& \sup_{k \in \RR} \int_M \varphi \, K_k \leq  \int_M \varphi \, \alpha_K. 
\endaligned
\ee

Next, we need to introduce a suitable generalization of Kruzkov and Kuznet\-sov's mollifiers. 
This is straightforward on a Riemannian manifold, by using the canonical distance function, but 
in the present formalism, we need the following new notion.  

In what follows, if $\zeta=\zeta(p,q)=\zeta_{p,q}$ is a function on
$M\times M$, then 
$d_p \zeta$ 
and $d_q \zeta$ 
denote its differentials with respect to the first and second arguments, respectively. 
We use the notation $\zeta_p$ for the function $p \mapsto \zeta(p,q)$, for $q$ fixed, and we 
often specify 
the integration variable in each integral to avoid confusion, by 
adding a subscript to the volume form under consideration. 
To a sequence of functions $\zeta^\delta$ we associate their supports $E_p^\param := \supp_q \zeta^\param_q$. 
Recalling that $\alphab$ denotes the volume form \eqref{DF.65}, 
we also write $|E| := \int_E \alphab$ for the volume of a set $E$
and, therefore,  
$|E_p^\param| = \int_{E_p^\param} \alphab_q$. We also use the notation $\dashint_E \alphab := |E|^{-1}\int_E \alphab$.

\begin{definition}
\label{DF-05} 
Fix a non-negative constant $A$  
 and a smooth $1$-form field $\beta$ defined on $M$. 
A family of non-negative (compactly supported and smooth) test-functions $(\zeta^\param)_{\param > 0}$ 
defined on $M\times M$ is called an {\bf $(A,\beta)$--admissible family of mollifiers}
if  the following conditions are
satisfied:
\enumerate 
\item {\bf Unit mass condition:}
\label{C.30}
 $\int_M \zeta_{p,q}^\param \, \bar\alpha_q = 1$ for each $p \in M$. 

\item {\bf Sup-norm condition:} 
 \label{C.10}  
$\sup_{q \in M} \zeta_{p,q}^\param \le |E_p^\param|^{-1}$ for each $p \in M$. 

\item {\bf Differential condition:}
\label{C.40} 
\[ \iint_{ M\times M}  |d_p\zeta_{p,q}^\param \wedge \gamma_{p,q}| 
\le 
\frac{1}{\param} \int_{M} \dashint_{E_p^\param} \, |\beta_p \wedge \gamma_{p,q}| 
\]
for each test $(2N+1)$-form field $\gamma$. 

\item {\bf Symmetry condition:}
\label{C.35} 
for each $\ubar \in [-C_0, C_0]$ and with $\gamma := \del_u \omega(\ubar)$ 
\[
\aligned
&\iint_{ M\times M} \varphi_{p,q} \, \Big(d_p \zeta_{p,q}^\param \wedge\gamma_p \wedge \alphabar_q + d_q 
\zeta_{p,q}^\param \wedge \gamma_q \wedge \alphabar_p \Big) 
\\
&
\leq  {A} \, \int_{M} \dashint_{E_p^\param} \varphi_{p,q} \, \alphabar_p \wedge 
\alphabar_q 
\endaligned
\]
for every bounded function $\varphi: M\times M \to \RR^+$.

\endenumerate 
\end{definition}

On a Riemannian manifold 
the distance function allows one to define such a family of test-functions. 
In general, the test-functions should be defined in each application by taking advantage of 
special properties of the given family of approximate solutions (cf.~examples below).

These assumptions arise from natural requirements on the supports $E_p^{\param}$, 
volumes $|E_p^\param|$, and test-functions $\zeta^\param$. They  
take the proposed form, due to the lack of metric structure on the manifold.
For instance, Condition \eqref{C.10} is a ``smallness" condition on the mollifiers, while Condition \eqref{C.30}
replaces the unit integral property of standard mollifiers in Euclidian or Riemannian manifolds.

Condition \eqref{C.35} is a symmetry property, enjoyed in the 
Euclidian space by  $\zeta^\param (p,q) =\eta(\ell(p,q))$, which takes then
 the much simpler form $d_p\zeta^\param = -d_q \zeta^\param$, 
when $\ell(p,q)$ denotes the Euclidian distance function.
In a Riemannian or Lorentzian setting, a similar 
property holds, but one needs to use parallel transport to compare $d_p \zeta$ and $d_q
\zeta$ (see the proof of Theorem \ref{10-10} below). Our condition above is intended to
encompass all situations, when no metric is naturally available on the manifold. 
The inequality in
Condition \eqref{C.35} in Definition~\ref{DF-05} is motivated by the following formal calculation, 
in wich $M=\RR$ and $\zeta$ coincides with a function $\eta$ composed with the distance function $\ell$:
$$
\aligned
d_p \zeta \gamma(p) + d_q \zeta \gamma(q) 
& = d_p \zeta (\gamma(p) - \gamma(q)) 
\\
& \le |\eta'(\ell(p,q))| \, \sup_\RR|\gamma'| \ell(p,q).
\endaligned
$$
Now, $\eta$ will usually be a standard mollifier, and thus in this example, $|\eta'| \le 
\param^{-2}\lesssim (\param |\supp \eta|)^{-1}.$ Since $\ell(p,q) \le \param$, the expression 
above is bounded by $|\supp \eta|^{-1} \sup_\RR |\gamma'|$, and 
Condition \eqref{C.35} is a generalization of this formal argument, in an integral form. 
This is 
also a condition on the size of the support of $\zeta^\param$, since it somehow generalizes 
the fact that in Riemannian space, $\ell(p,q) \le \param$ if $q$ is in the support of $\zeta^
\param_p$ and if those supports are geodesic balls of diameter $\param$. This assumption is 
necessary since, without a metric on the manifold, it seems difficult 
to reproduce the above argument in a geometric (i.e., coordinate independent) way.  

Note, however, that if one is interested in some particular problem, it is not hard to express 
the constant $A$ in Condition \eqref{C.35} as a (possibly coordinate dependent) quantity 
involving derivatives of $\gamma$.

As for Condition \eqref{C.40}
it amounts to a uniform upper bound on the form fields $d_p \zeta$ by a certain
1-form field $\beta$. It is analogous to the bound $|\nabla \zeta^\param| \le C \param^{-N-2}$ enjoyed by 
the standard molifiers in the Euclidian and Riemannian cases \cite{BP,LNO}.

%------------------------------------------------------------------------------------

\subsection{Statement of the error estimate}

Our main result in the present paper is now stated. Recall that all values $u$ under consideration belong 
to the 
interval $[-\boundd,\boundd]$.

\begin{theorem}[Error estimate for conservation laws on a spacetime]
\label{DF-10} 
Consider the conservation law \eqref{DF.10} with flux field $\omega$, 
posed on a spacetime $M$ satisfying the foliated condition \eqref{MHT} and the global hyperbolicity condition \eqref{GHC}
for some $\underline c, \overline c$.
Let $\zeta^\delta$ be an $(A,\beta)$--admissible family of mollifiers 
associated with some non-negative constant $A$ and $1$-form field $\beta$. 
Consider two functions 
$u(t), v(t): \Hcal_t \to \RR$ that belong to $L^1(\Hcal_t)$ for each $t \geq 0$ and 
are right-continuous in $t$.
Assume moreover that $v$ is an (exact) entropy solution to the conservation law \eqref{DF.60},
and that $u$ satisfies the approximate entropy inequalities \eqref{DF.40} for some $H_k, K_k, \alpha_H, \alpha_K$ satisfying the bounds 
\eqref{DF.50} and such that $i_{\Hcal_t}^*\alpha_H $ belongs to $L^1(\Hcal_t)$ 
for all $t$. 
Then, the following $L^1$-type estimate holds for all $\param >0$ and $T>0$
\be
\label{DF.100} \aligned
\int_{\Hcal_T} i^*\bfOmega (u, v) &\le \int_{\Hcal_0} i^*\bfOmega
(u, v)
+ R^\param[v] + R^\param[\omega] + R^\param[\alpha]  ,
\endaligned
\ee
where
$$
\aligned
R^\param[v] :=&  \sup_{t\in(0,T)} \int_{\Hcal_t} \, i^*
\alphab_p \,  
\dashint_{E_p^{\param}} 
|v_p -
v_q| \, B_q,
\endaligned
$$
$$
\aligned
R^\param[\omega]
:=& \int_{M_T} \dashint_{E_p^{\param}}  \, 
       \big|d\omega_p(v_q)\wedge \alphabar_q - d\omega_q(v_q)\wedge \alphabar_p \big|,
\\
R^\param[\alpha] 
:=& \frac{1}{\param} \int_{M_T} |\beta \wedge \alpha_H| +  \int_{\Hcal_0 \cup \Hcal_T} |i^* \alpha_H| + \int_{M_T} \alpha_K, 
\endaligned
$$ 
and $B$ is an $(N+1)$-form field in $q$ defined by
$$
B_q := \big(2 \, \cbar + T \, A \big) \, \alphabar_q + T \, \sup_u  | \del_u d\omega_q(u)|. 
$$
\end{theorem}

A few remarks about the above theorem are in order.
First of all, the terms $R^\param[v]$ and $R^\param[\omega]$
 are expected to tend to zero with $\param$. 
 For instance, when a metric is prescribed on the 
manifold, the term $R^\param[v]$ is estimated 
(see Lemma \ref{CL-10}, below) like in the classical Euclidian case~\cite{BP}:  
$R^\param[v] \le C \, T \, \param \TV(v(0))$, provided $v$ has bounded variation.  

Second, under the regularity assumptions on the flux $\omega$, the 
term $R^\param[\omega]$ is expected to be of order $\Ocal(\param)$. However, to
establish this property, 
one needs to control the ``size'' of the sets $E_p^{\param}$, 
but this cannot be formulated without a notion of distance on $M$. 
In contrast, Theorem \ref{10-10} below will specialize to the case of 
a metric on $M$ 
and on conservation laws based on vector fields, 
and we will see explicitly that $R^\param[\omega]$ vanishes with $\param$.

Finally, note that the quantity $\int_{\Hcal_T} i^*\bfOmega (u, v) $ can be seen as a measure of the $L^1$-norm of the 
difference between $u$ and $v$. Indeed, in the Euclidian and Riemannian cases, it reduces to 
$\int_{\Hcal_T} |u_p - v_p| d\mathrm{Vol}(\Hcal_t)$.

Before discussing some applications of the above theorem, it is interesting to consider 
the special case where the flux field $\omega$ is ``geometry compatible''. 

\begin{corollary}
\label{DF-20}
In addition to the assumptions in Theorem \ref{DF-10}, assume that the flux field $\omega$ is 
geometry compatible, in the sense that $(d\omega)(\ubar) =0$ for each $\ubar \in \RR$. Then, 
the following error estimate 
\be
\label{DF.102} \aligned
\int_{\Hcal_T} i^*\bfOmega (u, v) &\le \int_{\Hcal_0} i^*\bfOmega
(u, v)
+ R^\param[v]  + R^\param[\alpha],
\endaligned
\ee
holds for all $T \geq 0$, with $R^\param[v]$ and $R^\param[\alpha]$ defined as in Theorem~\ref{DF-10}. 

Furthermore, if $v$ is sufficiently smooth so that $R^\param[v] \to 0$ as $\param \to 0$ 
and if $u$ is also an exact entropy solution, then the $L^1$-like distance between $u$ and $v$ 
$$
t \mapsto \int_{\Hcal_t} i^*\bfOmega (u, v)
$$
is non-increasing in time. 
\end{corollary}

The second statement in the above corollary is nothing but
 the contraction property of the semi-group of entropy solutions. 

%------------------------------------------------------------------------------------------------------------------------------

\subsection{Application (I). Comparing two conservation laws}
\label{24s}
 
Theorem \ref{DF-10} applies to conservation laws with ``modified'' flux, and allows us to estimate the difference 
between entropy solutions to two distinct conservation laws. 
Let $\omega$ and $\widetilde\omega$ be two geometry-compatible 
flux fields, and introduce their corresponding Kruzkov entropy flux field 
$$
\Omega(v,k) = \sgn(v - k) (\omega(v) - \omega(k)), \qquad
\widetilde\Omega(u,k) = \sgn(u - k) (\widetilde\omega(u) - \widetilde\omega(k)). 
$$
The solutions $u,v$ under consideration satisfy the entropy inequalities
\be
\label{AP.10}
d( \Omega(v,k)) \le 0,\qquad 
d( \widetilde\Omega(u,k)) \le 0, \quad k\in \RR.
\ee
In order to avoid unnecessary technicalities, we may assume that the chain rule applies
to expressions involving the functions $u$ and $v$ (i.e.~bounded functions with bounded variation, for instance).

\begin{theorem}
\label{AP-10} Let $u,v$ be to entropy solutions satisfying \eqref{AP.10} for two flux fields
$\omega$ and $\widetilde \omega$, and 
assume that the conditions in Theorem~\ref{DF-10} hold for both conservation laws. Then,
 the following two estimates hold:
\begin{enumerate}
\item
If $v$ is sufficiently regular so that $R^\param[v] \to 0$ as $\param \to 0$, then 
\be
\label{AP.20}
\int_{\Hcal_{{T}}} i^*\bfOmega (u, v) \le \int_{\Hcal_{0}} i^*\bfOmega
(u, v) + C\int_{M_T} | \del_u (\omega - \widetilde\omega) \wedge dv|
\ee
for some uniform constant $C>0$. 

\item
If, moreover, $R^\param[v] \le \param \, \bar{R}[v]$, for some constant $\bar{R}[v]$ independent of $\param$, then
\be
\label{AP.30}
\aligned
\int_{\Hcal_{{T}}} i^*\bfOmega (u, v) 
\le 
& \int_{\Hcal_{0}} i^*\bfOmega
(u, v) + C \, \Bigg(  \bar{R}[v] \int_{M_T}|\beta \wedge Q(\omega,\widetilde\omega)| |u(p)| \Bigg)^{1/2}
\\
& + \int_{\Hcal_0 \cup \Hcal_T} i^*Q(\omega,\widetilde\omega) |u(p)|
\endaligned
\ee
\end{enumerate}
with
$$ 
{Q(\omega,\widetilde\omega) = \sup_{u\neq 0}|\omega(u) - \widetilde\omega(u)| / |u|}.
$$
\end{theorem}

\begin{proof} In order to apply Theorem~\ref{DF-10}, or more precisely its Corollary~\ref{DF-20}, 
we need to identify the structure of the relevant approximate conservation laws. To this end, we 
write 
\[
d(\Omega (u,k)) \le d\big( \sgn(u - k) \big( (\omega - \widetilde\omega)(u) -
(\omega - \widetilde\omega)(k) \big) \big) =: d(\gamma_k(u)).
\]
To show the estimate \eqref{AP.20}, we set $K_k := d(\gamma_k(u))$. Using a weak form of the chain rule, we
 see that
\[
|K_k| \leq |\del_u (\omega - \widetilde\omega) \wedge du |.
\] 
Hence, we arrive at the desired estimate \eqref{AP.20}, when $\param \to 0$ and 
after changing the role of $u$ and $v$. 

Second, to establish \eqref{AP.30}, we set $H_k := \gamma_k(u)$. 
Given an arbitrary 
$1$-form field and following arguments in \cite{BP} (say that $\omega(0) = \widetilde\omega(0) = 0$ for
 simplicity), we obtain 
\[
|H_k \wedge \gamma| \leq C \, |\alpha_H \wedge \gamma| 
\]
with $\alpha_H := Q(\omega,\widetilde\omega) \, |u|$. 
Therefore, the estimate \eqref{AP.30} now follows from Corollary \ref{DF-20} 
by minimizing over the parameter $\param$. This completes the proof of Theorem~\ref{AP-10}.
\end{proof}

%--------------------------------------------------------------------------------------------------------------- 

\subsection{Derivation of the error estimate}

\proof[Proof of Theorem \ref{DF-10}] Let $u$
satisfy the approximate entropy inequalities \eqref{DF.40}, and let
$v$ satisfy the (exact) entropy inequalities \eqref{DF.60}. Let
$\varphi=\varphi_{p,q}$ be a smooth, compactly supported function on
$M \times M$. According to \eqref{DF.20} and \eqref{DF.50}, for each $k\in \RR$ and $q \in M$ we have for $\varphi \geq 0$ 
\[
\aligned
& -\int_{M} d_p\varphi \wedge\bfOmega_p (u_p, k) + \int_{M} \varphi G_p(u_p, k) 
\\
&\le \int_{M} \varphi E_{k} = - \int_{M} H_k\wedge d_p\varphi + \int_{M} K_k \varphi
\\
& \le \int_{M} |\alpha_H \wedge d_p \varphi| + \int_{M} \alpha_K \varphi,
\endaligned
\]
thus by taking $k=v_q$ we find
\be
\label{DF.105}
\aligned
&-\int_{M} d_p\varphi \wedge\bfOmega_p (u_p, v_q) + \int_{M} \varphi G_p(u_p, v_q) 
\\ 
&\le  \int_{M} |\alpha_H \wedge d_p \varphi| + \int_{M} \alpha_K \varphi. 
\endaligned
\ee
On the other hand, taking $l=u_p$ in \eqref{DF.60} gives
\be
\label{DF.107}
\aligned
-\int_{M} d_q\varphi \wedge\bfOmega_q (u_p, v_q) + \int_{M} \varphi G_q(v_q,u_p) \le 0.
\endaligned
\ee
Since $q$ is the integration variable, the integrals in \eqref{DF.107} may be
viewed as real-valued functions of $p$. Therefore, we may integrate this
function on the manifold, provided a volume form is used. Likewise, 
we may integrate the inequality \eqref{DF.105} in $q$. Choosing
the form $\alphabar$ from \eqref{DF.65}, we obtain
\[
\aligned 
& -\int_M \alphabar_q \int_{M} d_p\varphi
\wedge\bfOmega_p (u_p, v_q) + \int_M \alphabar_q \int_{M}
\varphi G_p(u_p, v_q)
\\&
 \le \int_M \alphabar_q \int_{M} |\alpha_H \wedge d_p \varphi|
\endaligned
\]
and
\[
\aligned
-\int_M \alphabar_p\int_{M} d_q\varphi \wedge\bfOmega_q (u_p, v_q) + \int_M 
\alphabar_p\int_{M} \varphi G_q(v_q,u_p) \le 0.
\endaligned
\]
Summing up the above two inequalities and applying Fubini's theorem, we obtain
\be \label{DF.120} \aligned 
& -\iint_{M \times M} d_p\varphi
\wedge\bfOmega _p(u_p, v_q) \wedge \alphabar_q + d_q\varphi
\wedge\bfOmega _q(u_p, v_q) \wedge \alphabar_p
\\
& + \iint_{M \times M} \varphi (G_p(u_p, v_q) \wedge\alphabar_q +
G_q(v_q,u_p) \wedge\alphabar_p )
\\
&
\le \iint_{M \times M} |\alpha_H \wedge d_p \varphi|\wedge\alphabar_q +  \alpha_K \wedge\alphabar_q \varphi.
\endaligned
\ee

Let $(\zeta^\delta)_\param$ be an admissible family of mollifiers as in Definition \ref{DF-05}, and
let $\chi_p$ be a smooth,
compactly supported real function on $M$ to be specified later.
We choose the test-functions
$$
\varphi_{p,q} = \chi_p \zeta^\param_{p,q}, 
$$
which leads us to 
$$d_p \varphi = \zeta^\param d\chi + \chi d_p \zeta^\param, \qquad d_q \varphi = \chi d_q \zeta^
\param,$$
and so the inequality \eqref{DF.120} becomes
\be
\label{DF.140}
\aligned
& - \iint_{M \times M} \zeta^\param d \chi \wedge\bfOmega _p(u_p, v_q) \wedge \alphabar_q
\\
& - \iint_{M \times M} \chi \Big(  d_p\zeta^\param \wedge\bfOmega_p (u_p,
v_q) \wedge \alphabar_q + d_q\zeta^\param \wedge\bfOmega_q (u_p,
v_q) \wedge \alphabar_p \Big)
\\
& + \iint_{M \times M} \chi \zeta^\param (G_p(u_p, v_q) \wedge\alphabar_q +
 G_q(v_q,u_p) \wedge\alphabar_p )
\\
&
 \le \iint_{M \times M} |\alpha_H \wedge d_p(\chi\zeta^\param)|\wedge\alphabar_q + 
\alpha_K \wedge\alphabar_q \chi \zeta^\param,
\endaligned
\ee which, with obvious notation, has the form 
$$
I_1 -I_2 + I_3 \le I_4.
$$

\subsection*{The terms $I_2$ and $I_3$.}

We have 
\[
\aligned
I_2 &:=  \iint_{M \times M} \chi_p \, (  d_p\zeta^\param \wedge\bfOmega _p(u_p, v_q) \wedge 
\alphabar_q + d_q\zeta^\param \wedge\bfOmega _q(u_p, v_q) \wedge \alphabar_p)
\\
&= \iint_{M \times M} \chi_p \,  d_q \zeta^\param \wedge \bfOmega_q (u_p, v_p) \wedge \alphabar_p
\\
&\quad + \iint_{M \times M} \chi_p \, d_q \zeta^\param \wedge\Big( \bfOmega_q(u_p, v_q) - \bfOmega_q(u_p, v_p) \Big)\wedge \alphabar_p
\\
& \quad  + \iint_{M \times M} \chi_p \, d_p\zeta^\param \wedge \Big( \bfOmega_p(u_p, v_q)  - \bfOmega_p (u_p,v_p) \Big) \wedge \alphabar_q
\\
& =: I_{2,1} + I_{2,2} + I_{2,3},
\endaligned
\]
with obvious notations, where we have used 
\[
\aligned
& \iint_{M \times M}\chi_p d_p\zeta^\param \wedge  \bfOmega_p(u_p, v_p) \wedge \alphabar_q 
\\
&= \int_{M} \chi_p \bfOmega_p(u_p, v_p) \wedge d_p \Big( \int_M  \zeta^\param_q \alphabar_q \Big) = 0,
\endaligned
\]
by Condition \eqref{C.30} in Definition \ref{DF-05}.
We now analyze the terms in $I_2$. For the first term, note that
\[
\aligned
I_{2,1}&=  \int_{M} \chi_p \alphabar_p \Big( \int_{E_p^{\param}}  d_q \zeta^\param \wedge \bfOmega_q (u_p, 
v_p)  \Big) 
\\
&= -  \int_{M} \chi_p \alphabar_p \Big( \int_{E_p^{\param}}  \zeta^\param d \bfOmega_q (u_p, v_p)  
 \Big).
\endaligned
\]
This integration by parts is possible since $\bfOmega_q (u_p, v_p)$ depends on $p$ only 
through $u$ and $v$, i.e., the explicit spatial dependence is on $q$.

Next, we have
\[
\aligned
I_{2,2} =  \int_0^1 \iint_{M \times M} \chi_p \, (v_q - v_p) d_q \zeta^\param \wedge \del_v\bfOmega_q(u_p, v^*)  \wedge \alphabar_p \,ds,
\endaligned
\]
with $v^* = sv_q + (1-s) v_p$, and a similar expression for $I_{2,3}$. This gives
\[
\aligned
I_{2,2} + I_{2,3} &= \int_0^1 \iint_{M \times M} \chi (v_q - v_p)  \Big( d_q \zeta^\param \wedge  \del_v\bfOmega_q(u_p, v^*) \wedge \alphabar_p 
\\
&\qquad+ d_p \zeta^\param \wedge \del_v\bfOmega_p(u_p, v^*)\wedge \alphabar_q  \Big) ds,
\endaligned
\]
which, according to Condition \eqref{C.35} of Definition \ref{DF-05} and the definition of $\bfOmega$, yields, 
for some constant $A > 0$, 
\[
\aligned
I_{2,2} + I_{2,3} & \le \int_{M} \dashint_{E_p^\param} A \chi_p |v_q - v_p | \, \alphabar_p \wedge \alphabar_q.
\endaligned
\]
Putting these estimates together, we obtain
\be
\label{DF.180}
\aligned
I_2 
\le 
& -  \int_{M} \chi_p \alphabar_p \Big( \int_{E_p^{\param}}  \zeta^\param d \bfOmega_q (u_p, 
v_p)  \Big) 
\\
& + \int_{M} \dashint_{E_p^\param} A\chi_p |v_q - v_p|
\, \alphabar_p \wedge \alphabar_q.
\endaligned \ee

Now, we estimate the term $I_3$ in \eqref{DF.140}. First, 
note that
\[
\aligned
& G_p(u_p,v_q)\wedge\alphabar_q + G_q(v_q,u_p)\wedge\alphabar_p 
\\
&= \sgn(u_p - v_q) ( d\omega_p(v_q) \wedge \alphabar_q - 
d\omega_q(u_p) \wedge \alphabar_p ),
\endaligned
\]
and thus
\[
\aligned
& d\omega_p(v_q) \wedge \alphabar_q - d\omega_q(u_p) \wedge \alphabar_p  
\\
& = d\omega_p(v_q) \wedge \alphabar_q - d\omega_q(v_q) \wedge \alphabar_p
\\
&\quad+d\omega_q(v_q) \wedge \alphabar_p - d\omega_q(u_p) \wedge \alphabar_p
  =: A_1 + A_2.
\endaligned
\]
From Condition \eqref{C.10}, that is, $\zeta^\param \le |E_p^\param|^{-1}$ we find immediately 
\[
\aligned
\int_{M} \chi_p  \int_{E_p^{\param}} \zeta^\param \sgn(u_p - v_q)A_1  &\le \bar{R}^\param[\omega],
\endaligned
\]
where $\bar{R}^\param[\omega]$ is defined as $R^\param[\omega]$ in the statement of the theorem, but with ${\chi_p}/ |E_p^\param|$ instead of $1/|E_p^\param|$.
On the other hand, the term $A_2$ gives
\[
\aligned
\int_{M} \chi_p  \int_{E_p^{\param}} \zeta^\param \sgn(u_p - v_q) A_2  & = -\int_{M} \chi_p \alphabar_p 
\int_{E_p^{\param}} \zeta^\param d \bfOmega_q(u_p,v_q).
\endaligned
\]
This leaves us with
\[
\aligned
- I_3 &\le
 \int_{M} \chi_p \alphabar_p  \int_{E_p^{\param}} \zeta^\param  d\bfOmega_q(u_p,v_q)
 +\bar{R}^\param[\omega].
\endaligned
\]
From the last inequality and \eqref{DF.180}, we obtain
\[
\aligned
I_2 - I_3 &\le \int_{M} \chi_p\alphabar_p  \int_{E_p^{\param}} \zeta^\param \Big( d \bfOmega_q
(u_p,v_q) - d\bfOmega_q (u_p, v_p)\Big)
\\
& +\int_{M} \dashint_{E_p^\param} A \chi_p |v_q - v_p| \alphabar_p \wedge \alphabar_q + \bar{R}^\param[\omega].
\endaligned
\]
The first integral is bounded by
\[
\aligned \int_{M} \chi_p \alphabar_p  
\dashint_{E_p^{\param}} \sup_u | \del_u
d\omega_q(u)| |v_p - v_q| ,
\endaligned
\]
so that we finally find
\be
\label{DF.190}
\aligned 
&I_2  - I_3
\\
&{\le } \int_{M} A\alphabar_p \chi_p \dashint_{E_p^{\param}} 
\big( 
\alphabar_q + \sup_u | \del_u d\omega_q(u)| \big) |v_p - v_q|
\\&\qquad + \bar{R}^\param[\omega].
\endaligned
\ee 

\subsection*{The terms $I_1$ and $I_4$.}
Let us now turn to $I_1$, the main term in \eqref{DF.140}. As in \cite{BP}, we 
choose the test-functions $\chi =
\chi^\eps$ to be supported on $\bigcup_{0\le t\le T+\eps}
\Hcal_t$ and constant within any hypersurface $\Hcal_t$, so that
$\chi^\eps$ is a function of $t$ only. Moreover, we arrange that
$\chi^\eps\le 1$, $\chi^\eps \equiv 1$ on $\cup_{\eps \le t \le
T}\Hcal_t$, and that $\del_t\chi^\eps(t) \to \delta_{t=0} -
\delta_{t=T}$ as $\eps\to 0$, where $\delta_{t=\tau}$ is a Dirac mass centered at $\tau$. 
Also, we have $d\chi^\eps = \del_t \chi^\eps dt$.

We now find 
\[
\aligned
I_1 &= -\iint_{M \times M} \zeta^\param d \chi^\eps \wedge\bfOmega _p(u_p, v_q) \wedge \alphabar_q
\\
& = -\iint_{M \times M} \zeta^\param d \chi^\eps \wedge\bfOmega_p(u_p, v_p) \wedge \alphabar_q
\\
&\quad+ \iint_{M \times M} \zeta^\param d \chi^\eps \wedge \big( \bfOmega_p (u_p, v_p) - 
\bfOmega_p (u_p, v_q) \big)
\wedge \alphabar_q.
\endaligned
\]
Consider the last integral, and observe that from the non-degeneracy condition \eqref{GHC}
we find, for any positive real function $g$ supported in $M_T$,
\[
\aligned
&\int_M g_p |dt \wedge \del_u \bfOmega|  \le 2 \sup_u \int_0^T \int_{\Hcal_t} g_p i^*\del_u \omega(u)\,dt
\\
&\le 2\bar{c} \int_{M}g_p dt \wedge \del_u 
\omega(0)
= 2\bar{c}\int_M g_p \alphabar_p.
\endaligned
\]
Thus, using Condition \eqref{C.10},
\[
\aligned &\iint_{M \times M} \zeta^\param d \chi^\eps \wedge \big( \bfOmega_p (u_p,
v_p) - \bfOmega_p(u_p, v_q) \big) \wedge \alphabar_q
\\
&\qquad\le \int_M | d \chi^\eps \wedge \sup_v \del_v \bfOmega_p|
\dashint_{E_p^{\param}}  |v_p - v_q|
\alphabar_q
\\
&\qquad \le \int_M 2 | \del_t \chi^\eps| \bar{c} \, dt \wedge \del_u \omega_p(0) 
\dashint_{E_p^{\param}}  |v_p - v_q| \alphabar_q
\\
& \qquad=\int_M 2 | \del_t \chi^\eps| \bar{c}\, \alphabar_p 
\dashint_{E_p^\param}  |v_p - v_q| \alphabar_q.
\endaligned
\]
Furthermore, we have from Condition \eqref{C.30} in Definition \ref{DF-05},
\[
\aligned
-\iint_{M \times M} \zeta^\param d \chi^\eps \wedge\bfOmega_p (u_p, v_p) \wedge \alphabar_q = -\int_{M}  d 
\chi^\eps \wedge\bfOmega_p (u_p, v_p).
\endaligned
\]
Therefore, recalling \eqref{DF.140} and \eqref{DF.190}, we find
\be
\label{DF.200}
\aligned
&-\int_{M}  d \chi^\eps \wedge\bfOmega_p (u_p, v_p)
\\
&\quad \le \int_M 2 | \del_t \chi^\eps| \bar{c}\, \alphabar_p 
\dashint_{E_p^{\param}}  |v_p - v_q| \alphabar_q
\\
& \qquad + \int_{M}  A \alphabar_p \chi^\eps \dashint_{E_p^{\param}} 
\big( \alphabar_q + \sup_u | \del_u d\omega_q| \big) |v_p - v_q| 
\\
&\qquad + \bar{R}^\param[\omega] + I_4.
\endaligned
\ee

\medskip

Let us now take $\eps\to 0$. First, we have
\[
\aligned
&-\limsup_{\eps \to 0}\int_{M}  d \chi^\eps \wedge\bfOmega_p (u_p, v_p) 
\\
&\qquad = - \limsup_{\eps \to 0}
\int_0^{T+
\eps}  \del_t \chi^\eps(t) \int_{\Hcal_t} i^*\bfOmega_p (u_p, v_p) dt
\\
&\qquad \le\int_{\Hcal_T} i^*\bfOmega_p (u_p, v_p) - \int_{\Hcal_0} i^*\bfOmega_p (u_p, 
v_p).
\endaligned
\]
Similarly, and since $|\del_t \chi^\eps(t)| \to \delta_{t=0} + \delta_{t=T}$ as $\eps \to 0$, 
we find
\[
\aligned
&\limsup_{\eps \to 0} \int_M 2 | \del_t \chi^\eps| \bar{c}\, \alphabar_p 
\dashint_{E_p^{\param}}  |v_p - v_q| \alphabar_q
\\
& \qquad + \int_{M} A \alphabar_p \chi^\eps  \dashint_{E_p^{\param}} 
\big(  \alphabar_q + \sup_u | \del_u d\omega_q| \big) |v_p - v_q|
\\
& \le R^\param[v],
\endaligned
\]
with $R^\param[v]$ defined as in the statement of the theorem. Also, it is clear that  $\limsup_{\eps \to 0} \bar{R}^\param[\omega] \le R^\param[\omega]$.

 Finally, we must deal with the term
\[
I_4 = \iint_{M \times M} |\alpha_H \wedge d_p(\chi^\eps\zeta^\param)|\wedge\alphabar_q + 
\alpha_K \wedge\alphabar_q \chi^\eps \zeta^\param.
\]
We have, using Condition \eqref{C.40} in  
Definition \ref{DF-05} and \eqref{DF.50},
\[
\aligned
I_4 & 
 \le  \int_M |d\chi^\eps \wedge \alpha_H | \int_M \zeta^\param \alphabar_q +
\frac1\param \int_{M} {\chi^\eps} |\beta_p \wedge \alpha_H| 
\dashint_{E_p^{\param}} \alphabar_q
\\
& + \int_M \chi^\eps \alpha_K \int_M \zeta^\param \alphabar_q.
\endaligned
\]
Using Condition \eqref{C.30} and taking the $\limsup$ as $\eps \to 0$ gives $I_4 \le R^\param[\alpha]$. 
This completes the proof of Theorem \ref{DF-10}.
\endproof 

%======================================================

\section{Error estimates for a Lo\-rentzian manifold}
\label{LM}

\subsection{Conservation laws based on vector fields}

In this section, we derive error estimates for conservation laws posed on Lorentzian manifolds. 
This is motivated by the fact that Theorem~\ref{DF-10} is greatly improved
if the manifold under consideration has a metric defined on it. To begin with, one does not need to
assume the existence of a special family of mollifiers, since these are naturally provided by the metric.
Second, one can introduce second order error terms which allow us to 
consider more general approximate solutions to conservation laws; namely, we obtain
an error estimate for very general nonlinear diffusion models; see Theorem~\ref{ND-10} below.

Let $(\Lorentz,g)$ be a time-oriented, $(N+1)$-dimensional
Lorentzian manifold. Here, $g$ is a metric with signature $(-, +,
\ldots, +)$, and we recall that tangent vectors $X \in
T_p\Lorentz$ at a point $p \in \Lorentz$ can be separated into
timelike vectors ($g(X,X) < 0$), null vectors ($g(X,X) = 0$), and
spacelike vectors ($g(X,X) > 0$). The manifold is assumed to be
time-oriented, so that we can distinguish between past-oriented
and future-oriented vectors. The Levi-Civita connection associated
to $g$ is denoted by $\nabla$ and, for instance, allows us to
define the divergence operator $\dive$. Finally, we denote by
$dV_g$ (or $dV_g(p)$, to stress the integration variable) the volume
measure determined by the metric $g$.

Following \cite{BL}, a \emph{flux-vector} on a manifold is defined as a vector field 
$f=f_p(\ubar)$ depending on a real parameter
 $\ubar$ and the \emph{conservation law} on $(\Lorentz,g$) associated with $f$ reads
\be
\dive \big( f_p(u) \big) = 0, \qquad u: \Lorentz \to \RR.
\label{CL.1}
\ee
Moreover, the flux-vector $f$ is said to be \emph{geometry compatible} if
\be
\dive f_p (\ubar) = 0, \quad \ubar \in \RR, \, p \in \Lorentz,
\label{CL.2}
\ee
and to be \emph{timelike} if its $u$-derivative is a timelike vector field
\be
g\big(\del_u f_p(\ubar), \del_u f_p(\ubar) \big) < 0, \qquad p \in \Lorentz, \, \ubar \in 
\RR.
\label{CL.3}
\ee

We are interested in the initial-value problem associated with \eqref{CL.1}. So, we fix a 
spacelike hypersurface
$\Hcal_0 \subset \Lorentz$ and a measurable and bounded function $u_0$ defined on $
\Hcal_0$.
Then, we search for a function $u=u_p \in L^\infty(\Lorentz)$ satisfying \eqref{CL.1}  in the 
distributional
sense and such that the 
trace of $u$ on $\Hcal_0$ coincides with $u_0$, that is \be
u_{|\Hcal_0} = u_0. \label{CL.4} \ee It is natural to require that
the vectors $\del_u f_p(\ubar)$, which determine the propagation
of waves in solutions of  \eqref{CL.1}, are timelike and
future-oriented. Thus, we will assume throughout that the
flux-vector in equation \eqref{CL.1} is timelike, in the sense of
\eqref{CL.3}.

As in the previous section, we assume that the manifold $\Lorentz$ is \emph{globally
hyperbolic,} which in this Lorentzian setting means that there exists a foliation of
$\Lorentz$ by spacelike, compact, oriented Riemannian
hypersurfaces $\Hcal_t$ ($t \in \RR$):
$$
\Lorentz = \bigcup_{t \in \RR} \Hcal_t.
$$
Any hypersurface $\Hcal_{t_0}$ is referred to as a \emph{Cauchy surface} in $\Lorentz$,
while the family $\Hcal_t$ ($t \in \RR$) is called an \emph{admissible foliation associated 
with} $\Hcal_{t_0}$.
The future of the given hypersurface will be denoted by
$$
\Lorentzplus := \bigcup_{t \geq 0} \Hcal_t.
$$
Moreover, we 
denote by $n^t$ the future-oriented, normal vector
field to each $\Hcal_t$, and by $g^t$ the induced metric. Finally,
along $\Hcal_t$, we denote by $X^t$ the normal component of a
vector field $X$, thus $X^t := g(X,n^t)$. In the following, when
there is no risk of confusion, we write $F(u)$ instead of
$F_p(u)$.

\begin{definition}
\label{CL-0} A flux $F=F_p(\ubar)$ is called a \emph{convex
entropy flux} associated with the conservation law \eqref{CL.1} if
there exists a convex function $U:\RR \to \RR$ such that
$$
F_p(\ubar) = \int_0^\ubar \del_u U(u') \,\del_u f_p(u') \, du', \qquad p \in \Lorentz, \, 
\ubar \in \RR.
$$
A measurable and bounded function $u=u_p$ is called an \emph{entropy solution}
of the conservation law \eqref{CL.1}--\eqref{CL.2}
if the following \emph{entropy inequality}
$$
\aligned
& \int_{\Lorentzplus} g(F(u), \nabla_{g} \phi) \, dV_g + \int_{\Lorentzplus} (\dive F)(u) \, 
\phi \, dV_g
\\
& + \int_{\Hcal_0} g_0(F(u_0), n^0) \, \phi_{\Hcal_0} \, dV_{g_0}  -\int_{\Lorentzplus} 
U'(u) (\dive f)(u)\, \phi\, dV_{g} \geq 0
\endaligned
$$
holds for all convex entropy flux $F=F_p(\ubar)$ and all smooth functions $\phi \geq 0$ 
compactly supported in
$\Lorentzplus$.
\end{definition}

In particular, the requirements in the above definition imply the inequality
\be
\label{CL.10}
\dive \big( F(u) \big) - (\dive F)(u)  + U'(u) (\dive f)(u) \leq 0
\ee
in the distributional sense. For well-posedness results for the initial value problem
\eqref{CL.1}--\eqref{CL.4}, see \cite{BL,LO}.

%--------------------------------------------------------------------------------------------------------------

\subsection{Statement of the error estimate}

For convenience, we consider a Riemannian metric $\gbar$
associated with  the Lor\-entz\-ian metric $g$. We fix the natural
one, that is, in local coordinates where the matrix of the
metric $g$ is diagonal, we set $\gbar_{11}:= - g_{11}$ and
$\gbar_{ii}:= g_{ii}$, $(i = 2, \dots, N+1)$. For instance, this allows us 
to consider the distance function $\ell_\gbar$ associated
with $\gbar$. In particular, the volume form and divergence operator associated to $g$ or $\gbar$
are the same. Also, we write $B_p(r)$ for the geodesic ball centered at  $p\in M$ with radius $r$,
with respect to the metric $\gb$.

Since we shall rely on Kruzkov's family of entropies for the statement of our results as well as 
the proofs, we write the conservation law \eqref{CL.10} with Kruzkov's entropy flux
\be
\label{10.05}
F_p(u,k) = \sgn(u-k) (f_p(u) - f_p(k)), \qquad k\in \RR.
\ee

Thus, we are given a bounded measurable function $u$ satisfying an \textsl{approximate 
entropy inequality},
\be
\label{10.10}
\aligned
&\dive \big( \sgn(u_p  - k)( f_p(u_p) - f_p(k)) \big) + \sgn(u_p  - k) \dive f_p(k)
\\
& \qquad\le \dive H_k + K_k + \divex_{g_t}( a_k \nabla_{g_t} L_k).
\endaligned
\ee 
Here, the error terms $H_k, K_k, L_k$ are defined as follows: for each $k\in \RR$,
$H_k$ is a distributional vector field, that is, an element of the 
space of linear functionals from the space of smooth 1-forms and taking values in the space of 
(scalar) distributions on $M$. Thus, 
for each $\gamma \in \Lambda^1(M)$, $\langle H_k , \gamma \rangle \equiv \gamma(H_k)$ is 
a distribution on $M$, which we assume to be a Radon measure.

The terms $K_k$ and $L_k$ are Radon measures, and $a_k$ are 
continuously differentiable functions defined 
on $M$. We suppose that $H_k$, $K_k$, and $L_k$ satisfy the following uniform bounds (with respect 
to $k$),  
\be 
\label{10.15} 
|\gamma(H_k)| \le \alpha_H | \gamma^\sharp|,
\qquad
|K_k| \le \alpha_K, 
\qquad
|L_k| \le \alpha_L,
\ee
for some positive Radon measures $\alpha_H$, $\alpha_K$, $\alpha_L$ on $M$. Here, if $\alpha$ 
is a measure, $|\alpha|$ 
denotes its variation in the measure-theoretic sense, and $\gamma^\sharp$ 
is the vector obtained from the 1-form $\gamma$ by raising indices using the metric.
We also assume that the functions $a_k$ satisfy, for some $ \alpha_a$ independent of $k$,
\[
|a_k|, |\nabla_{g^t} a_k|_{g^t} \le \alpha_a.
\]

Note that due to the presence of a volume form, a measure on $M$ may be seen indifferently as an 
element of the dual of the space of $(N+1)$-forms or of the dual of the space of test-functions on 
$M$. Thus, when we write, for instance, $\int_M \alpha_K \varphi dV_{\gb}$, this denotes the duality 
between the (scalar) distribution $\alpha_K$ and the $(N+1)$-form $\varphi dV_{\gb}$, and so no 
regularity is assumed on $\alpha_K$.

Now, let $v$ denote the exact solution to the conservation law
\eqref{CL.10}, i.e. for all $l \in \RR$ 
\be \label{10.20} \aligned &\dive \big( \sgn(v_q - l)( f_q(v_q)
- f_q(l) ) \big) + \sgn(v_q  - l) \dive f_q(l) \le 0.
\endaligned
\ee
Defining
\[
G_p(u,k) = \sgn(u-k) \dive f_p (k),
\]
the entropy inequalities read
\be
\label{10.30}
\dive (F_p(u_p,k) ) + G_p(u_p,k) \le \dive H_k + K_k + \divex_{g_t}( a_k \nabla_{g_t} L_k)
\ee
for the approximate solution $u$, and
\be
\label{10.40}
\dive (F_q(v_q,l) ) + G_q(v_q,l) \le 0
\ee
for the exact solution $v$.

Our main result in this section is Theorem \ref{10-10}. It gives a precise, quantitative estimate 
of the evolution of
\[
\int_{\Hcal_t} F^t(u_p, v_p) \, dV_{g^t}.
\]
Note that this quantity is equal to
\[
\int_{\Hcal_t} |u_p - v_p| dV_{g^t}
\]
whenever the flux function $f$ of the equation is such that
$f^t(u) = u$ for all $u$. Therefore our estimates have the same
form as the usual estimates in \cite{BP}, where the manifold is
flat and the time evolution trivial. In the general case, we have
$F^t(u,v) = |f^t(u) - f^t(v)|$ for every $u,v\in \RR$. Since
$f^t$ is, by assumption, a strictly monotone function of $u$, this
quantity provides an equivalent measure of the difference between
$u$ and $v$ in the $L^1$-norm, which takes into account the geometry of the manifold
and the structure of the time-evolution of the foliation under
consideration.

\begin{theorem}[Error estimate for conservation laws on a Lorentzian manifold]
Let $u$ be a function satisfying 
the \label{10-10} approximate entropy inequalities \eqref{10.10}, \eqref{10.15}, and let 
$v$  be an 
exact solution satisfying \eqref{10.20}. Suppose also that
$u,v$ are right-continuous with values in $L^1 (\Hcal_t)$
and that,
for some $T\ge0,$
$\alpha_H$ is right-continuous from
$[0,T)$ with values in $L^1 (\Hcal_t)$. For
$u\in\RR$ and $p\in M$, define the constants
 \be \label{10.41}
\aligned
& \Lambda_0 := \sup_p \Lip_u f^t,
&{}
&\Lambda_1 := \sup_p \sup_{{X\in T_pM}\atop{|X| =1}} \Lip_u(\nabla_X f),
&
\\
&\Lambda_2 := \sup_u \Lip_p(\dive f), 
&
&\Lambda_3 := \sup_p \Lip_u(\dive f). 
\endaligned  
\ee
Then, for all $\param>0$ the estimate
\be \label{10.42}
 \aligned 
& \int_{\Hcal_T} F^t(u_p, v_p) \,
dV_{g^T} 
\\
& \le  \int_{\Hcal_0} F^t(u_p, v_p)
\, dV_{g^0} + C(E^\delta_v+ E^\delta_f  + E^\delta_H + E^\delta_K + E^\delta_L)
\endaligned
\ee
holds, 
where $C$ is a constant (which may depend on $N$) and
\[
E^\delta_v := (T\Lambda_1 + T\Lambda_3 + \Lambda_0)  \sup_{t \in (0,T)}  \int_{\Hcal_t} 
\dashint_{B_p(\param)}
|v_p- v_q| \,dV_g\,dV_{g^t},
\]
\[
E^\delta_f := T\sup_{t\in (0,T)} |\Hcal_t| \param\Lambda_2, 
\quad 
E^\delta_H := \int_{\Hcal_0 \cup \Hcal_T} \alpha_H dV_{g^t} +
\frac{1}{\param} \int_0^T\int_{\Hcal_t} \alpha_H dV_{g^t} dt,
\]
\[
E^\delta_K :=\int_0^T\int_{\Hcal_t} \alpha_K dV_{g^t} dt,
\qquad 
E^\delta_L :=  \frac{1}{\param^2} \int_0^T\int_{\Hcal_t} \alpha_L\alpha_a dV_{g^t} dt.
\]
\end{theorem}

\subsection{Application (II). The semi-group contraction property} 
\label{33s}

Our main result implies a key property of the semi-group of entropy solutions. 

\begin{corollary}\label{10-11}
In addition to the assumptions of Theorem \ref{10-10}, suppose also that the flux $f$ is geometry-compatible and that $v$ has bounded total variation on each slice $\Hcal_t$. Then, for 
every $\param >0$, the estimate
\[
\aligned 
& \int_{\Hcal_T} F^t(u_p, v_p) dV_{g^T} \,
\\
&  \le 
 \int_{\Hcal_0} F^t(u_p, v_p) dV_{g^0}
+ C T \, \param \, \sup_t \TV(v_{|\Hcal_t}) + C \, (E^\delta_H + E^\delta_K + E^\delta_L) 
\endaligned
\]
holds where $C$ is independent of $\param$, but depends on $f$, and the error terms are 
defined as in Theo\-rem 
\ref{10-10}. In particular, if $u$ is also an exact solution, the function
\[ 
t \mapsto \int_{\Hcal_t} F^t(u_p, v_p) \, dV_{g^t}
\]
is non-increasing.
\end{corollary}

This result is an immediate consequence of Theorem \ref{10-10} and of the 
following result, which provides the link between the term $E^\delta_v$ (associated with the 
regularity 
of the exact solution $v$) and its total variation. 
Recall that for each $p\in M$, the \emph{exponential map}
$
\exp_p : B_0(\delta) \to B_p (\delta) \subset M
$
provides a diffeomorphism between the ball of radius $\delta$ on the tangent space at $p$, 
$B_0(\delta) \subset T_pM$, and the geodesic ball (according to the Riemannian metric $
\gb$) $B_p(\delta)$ around $p$. 
(Here, $\delta$ must be small enough.)  
The exponential map provides a local chart by identifying $T_pM$ with $\RR^{N+1}.$ 
Also, in what follows, we write $dp,$ $dq$ 
instead of 
the volume element on $M$ to keep the exposition uncluttered and to stress the integration 
variable. It is also more convenient to write $v(p)$ instead of $v_p$.

\begin{lemma}
\label{CL-10}
Let $v \in BV(\Hcal_t), t\ge 0$ be a solution to the conservation law 
\eqref{10.20}. For $\param$  
sufficiently small, one has 
\[
\aligned
& \int_{\Hcal_t} \dashint_{B_p(\param)} |v(p)  - v(q)| \,dq\,dp 
\\
&\le
C \param \, 
\Big( \big(1 + \frac{ \| \Lip_u f\|_{L^\infty(M)}}{\beta}\big) \, \TV(v_{|\Hcal_t}) +  \frac{1}{\beta}\| (\dive f)_{|\Hcal_t} \|
_{L^1(\Hcal_t)}  \Big),
\endaligned
\]
where $\beta = \inf_{u,p}\del_u f^t (u,p) > 0$. In particular, if the flux $f$ is 
geometry-compatible, then
\[
\aligned
\int_{\Hcal_t} \dashint_{B_p(\param)} |v(p)- v(q)| \,dq\,dp &\le
C \param(1 + \frac{ \| \Lip_u f \|_{L^\infty(M)}}{\beta})   \TV(v_{|\Hcal_t}).
\endaligned
\]
\end{lemma}

\begin{proof}
We may assume that the function $v$ is sufficiently smooth since, 
by a 
standard density argument, 
the general result then follows for all functions $v$ with bounded variation.
Using the exponential map, we write 
\[
\aligned
\int_{\Hcal_t} \dashint_{B_p(\param)} |v(p)- v(q)| \,dq\,dp &=  \int_{\Hcal_t} 
\dashint_{B_0(\param)} |v(\exp_p(0))- v(\exp_p(h))| \,dh\,dp.
\endaligned
\]
Now, consider a partition of unity $\psi_i$, $1\le i \le m$ subordinate to a covering $\tilde 
U_i$ containing a neighborhood of $\Hcal_t$ of radius $\param$, and write $U_i = \Hcal_t 
\cap \tilde U_i$. We have
\[
\aligned
&\int_{\Hcal_t} \dashint_{B_0(\param)} |v(\exp_p(0))-  v(\exp_p(h))| \,dh\,dp
\\
&= \sum_i \int_{U_i} \psi_i \dashint_{B_0(\param)} |v(\exp_p(0))- v(\exp_p(h))| \,dh\,dp
\\
& \le  \dashint_{B_0(\param)}\sum_i \int_{U_i} \psi_i |v(\exp_p(0))- v(\exp_p(h))| \,dp\,dh
\\
&\le C(M) \sup_{h\in B_0(\param)}\sum_i \int_{U_i} \psi_i  |v(\exp_p(0))- v(\exp_p(h))| \,dp,
\endaligned
\]
where $C(M)$ is a constant depending on the geometry of $M$. Note that it is not trivial to reverse the order of integration above, and it is necessary to use 
the partition of unity.
This is due to the fact that there may be no way to globally specify the isomorphisms 
between $B_0(\param) \subset \RR^{N+1}$ and $B_0(\param) \subset TM$ in a smooth way. 
Indeed, a sufficient condition to be able to do so is that the tangent bundle of $M$ is trivial 
(i.e., diffeomorphic to $M\times TM$) in a neighborhood of $\Hcal_t$, which may not be 
the case.
Also, $\param$ must be small enough so that for each $p$ the 
point $\exp_p(h)$ is well defined.
We have also abused the notation somewhat since ``$dh$" really  stands for the determinant of the Jacobian of the exponential map, which, by compactness,  may be uniformly bounded as a function of $p$, whence the constant $C(M)$.
 
 Thus, using that $\exp_p(0) = p$, we find
\[
\aligned
\sum_i \int_{U_i} \psi_i |v(\exp_p(0))- v(\exp_p(h))| \,dp &\le \sum_i \int_{U_i} \psi_i  
\int_0^1| \frac{d}{ds} v(\exp_p(sh)| \,ds\,dp
\\
&\le  \param \int_0^1 \sum_i \int_{U_i} \psi_i  | \nabla_{\gb} v(\exp_p(sh))| \,dp\,ds
\\
&\le  \param  \int_{\Hcal_t} | \nabla_{\gb} v(p)| \,dp + A,
\endaligned
\]
with
\[
A =  \param \int_0^1 \sum_i \int_{U_i} \psi_i  | \nabla_{\gb} v(\exp_p(sh)) - \nabla_{\gb} 
v(\exp_p(0)) | \,dp \, ds.
\]
Next, split the gradient of $v$ into its time and spatial components,
\[
\aligned
\int_{\Hcal_t} | \nabla_{\gb} v(p)| \,dp &\le  \int_{\Hcal_t} | \del_t v(p)| \,dp + 
\int_{\Hcal_t} | \nabla^x_{\gb} v(p)| \,dp.
\endaligned
\]
Now we use the conservation law to estimate the temporal gradient in terms of the spatial 
gradient. Consider any system of local coordinates on the leaf $\Hcal_t$. A simple 
computation shows that the conservation law reads for smooth solutions
\[
\del_u f_p^t(v(p)) \del_t v(p) + \langle \del_u f_p(v(p)), \nabla_\gb^x v(p) \rangle_\gb + \dive f_p 
(v(p)) = 0.
\]
Thus, after integrating on $\Hcal_t$ we find 
\[
\aligned
\int_{\Hcal_t} |\del_t v| dp 
& \le \frac{1}{\beta} \int_{\Hcal_t} \del_u f^t |\del_t v| dp 
\\
& \le
\frac{ \| \Lip_u f \|_{L^\infty(M)}}{\beta} \int_{\Hcal_t} | \nabla^x_{\gb} v|dp + \frac{1}{\beta}\int_{\Hcal_t} |
\dive f_p (v(p)) | dp.
\endaligned
\]

Finally, since for smooth $v$
\[
\int_{\Hcal_t} | \nabla^x_{\gb} v(p)| \,dp = \TV(v_{|\Hcal_t}),
\]
the result will be proved if $\sup_{h \in B_0(\param)}A=\Ocal(\param^2)$, which is 
straightforward. This completes the proof of Lemma~\ref{CL-10}.
\end{proof}

%---------------------------------------------------------------------

\subsection{Application (III). A nonlinear diffusion model}
\label{ND}

In this section, we apply our results to a nonlinear diffusion model on a Lorentzian manifold.
For simplicity, we will consider the geometry-compatible case, in which the divergence of 
the
flux vanishes. Also, to shorten the presentation, we assume here that $u$ and $v$ are regular 
enough so that a weak form of the chain rule applies, see \cite{BP} where (in the Euclidian case) the 
required 
regularity is that $u,v \in BV_\loc$, which is also our case.
Following \cite{BP}, let $\phi : \RR \to \RR$ be a Lipschitz continuous function. 
We consider the nonlinear diffusion equation on $M$,
\[
\dive (f_p(u)) = \Delta_{g^t} \phi(u),
\]
where $\Delta_{g^t}$ denotes the Laplace-Beltrami operator on the leaf $\Hcal_t$.  One can check that if $\phi$ is non-decreasing, the approximate entropy inequalities are
\be
\label{ND.10}
\aligned
&\dive \big( \sgn(u_p  - k)( f_p(u_p) - f_p(k)) \big) 
\le \Delta_{g^t} |\phi(u) - \phi(k)|.
\endaligned
\ee
We will obtain, if $u$ and 
$v$ have bounded total variation, an error estimate 
in $\sqrt{\Lip \phi}$, which is the usual estimate in $\sqrt\eps$ when $
\phi(u) = \eps u$, and also a finer estimate in which only $v$ is required to have bounded total variation, and $u$ is only required to be bounded in $L^1$.

\begin{theorem}
\label{ND-10}
Let $T>0$ and assume 
 that $u$ satisfy the approximate entropy inequalities \eqref{ND.10} with a Lipschitz 
continuous nonlinear viscosity $\phi$. Let $v$ be an exact solution of the conservation law 
\eqref{10.20}. 
\begin{enumerate}
\item
Suppose that $\TV(v|_{\Hcal_t}) \le V$, $\TV(u|_{\Hcal_t}) \le 
U$ for all $0\le t\le T$. Then there is a constant $C>0$ such that
\be
\label{ND.20}
\aligned 
&\int_{\Hcal_T} F^t(u_p, v_p) \,
dV_{g^T} 
\le \int_{\Hcal_0} F^t(u_p, v_p)
\, dV_{g^0} + CT \sqrt{(\Lip \phi) \, V U }.
\endaligned
\ee
\item
Suppose that $\TV(v|_{\Hcal_t}) \le V$, and that $\int_{M_T} |u| dV_g \le U.$ Then,
\be
\label{ND.30}
\aligned 
&\int_{\Hcal_T} F^t(u_p, v_p) \,
dV_{g^T} \le \int_{\Hcal_0} F^t(u_p, v_p)
\, dV_{g^0} + CT \,  ( Q U)^{1/3} V^{2/3},
\endaligned
\ee
with $Q = \sup_{u\neq 0} |\phi(u) - \phi(0)| / |u|$.
\end{enumerate}

\end{theorem}

\begin{proof}
Apply Corollary \ref{10-11} with $H_k$ as the sole error term. We find (in the sense of measures)
\[
|H_k | \le \Lip_u \phi |\nabla_g u| =: \alpha_H.
\]
Estimating the term $R^\param[v]$ using Lemma \ref{CL-10} gives
\[
\aligned
& \int_{\Hcal_T} F^t(u_p, v_p) \,
dV_{g^T} 
\\
& \le \int_{\Hcal_0} F^t(u_p, v_p)
\, dV_{g^0} 
+ CT\param \, V + C \frac{T}\param \Lip \phi\, U,
\endaligned
\]
and the estimate \eqref{ND.20} follows by minimizing with respect to $\param$.

To establish \eqref{ND.30}, write $L_k = |\phi(u) - \phi(k)| - |\phi(0) - \phi(k)|$, and so 
$|L_k| \le Q |u| =: \alpha_L$. Now we apply Corollary \ref{10-11} and Lemma \ref{CL-10}  
to find
\[
\aligned
& \int_{\Hcal_T} F^t(u_p, v_p) \,
dV_{g^T} 
\\
& \le \int_{\Hcal_0} F^t(u_p, v_p)
\, dV_{g^0} 
+ CT\param \, V + C \frac{T}{\param^2} Q\, U,
\endaligned
\]
from which the estimate \eqref{ND.30} follows by choosing the optimal value of $\param$.
This completes the proof of Theorem \ref{ND-10}.
\end{proof}

%---------------------------------------------------------------------------------------------------------

\subsection{Derivation of the error estimate} 
\label{230} 

We provide here the proof of Theorem~\ref{10-10}. 
Let $\varphi$ be a smooth, compactly supported function on $M \times M$. From \eqref{10.15},
\eqref{10.30} and \eqref{10.40}, we have for all $k,l \in\RR$
\[
\aligned
&-\int_M d_p \varphi \big( F_p (u_p,k)\big) dV_g(p)
+ \int_M \varphi_{p,q} G_p(u_p,k) dV_g(p)
\\
& \le  \int_M (\dive H_k   +   K_k   + 
 \divex_{g^t}(a_k \nabla_{g^t} L_k) )\varphi_{p,q}  \,dV_g(p) 
\\
&\le  \int_M |d\varphi ( H_k) | \,dV_g(p) + \int_M |K_k| \varphi\,dV_g(p) + \int_M |L_k| \divex_{g^t}
(a_k \nabla_{g^t} \varphi) \,dV_g(p)
\\
&\le \int_M E(\varphi) \,dV_g(p),
\endaligned
\]
with
\be
\label{10.43}
\aligned
E(\varphi) :=
\alpha_H |\nabla_g \varphi|  + \alpha_K \varphi 
+  \alpha_L \alpha_a (  |\Delta_{g^t} \varphi| + |\nabla_{g^t} \varphi|),
\endaligned
\ee
and similarly
\[
\aligned &-\int_M d_q \varphi \big( F_q (v_q,l) \big) dV_g(q) + \int_M
\varphi_{p,q} G_q(v_q,l) dV_g(q) \le 0,
\endaligned
\]
where we recall that $d_p$ (resp.,~$d_q$) denotes the differential of a
function only with respect to $p$ (resp., $q$).
Taking $k=v_q$, $l=u_p$, integrating over $M$, summing the
above inequalities, and using Fubini's theorem, we find
\be \label{10.45} \aligned 
& - \iint_{M \times M}  \Big\{ d_p \varphi
\big(F_p(u_p,v_q) \big) + d_q \varphi\big( F_q(u_p,v_q)\big) \Big\}
\, dV_g(p)\,dV_g(q)
\\
&\quad + \iint_{M \times M} \varphi_{p,q} \big( G_p(u_p,v_q) + G_q(v_q,u_p) \big) \,dV_g(p)\,dV_g(q)
\\
& \le \iint_{M \times M} E(\varphi) \,dV_g(p)\,dV_g(q). 
\endaligned
\ee
Note also that
\[
\aligned
& G_p(u_p,v_q) + G_q(v_q,u_p)
\\
& = \sgn(u_p - v_q) ( \dive f_p(v_q) - \dive f_q(u_p)).
\endaligned
\]

Now take
\[
\varphi_{p,q} = \chi_p \; \xi^\param_{p,q},
\]
where the function $\xi^\param$ is of the form $\xi^\param_{p,q}= \zeta^\param(
\ell_\gb(p,q))$, the function $\chi$ is to be chosen later, 
and $\ell_\gb$ denotes the distance function associated with the Riemannian metric $\gb$.
Also, we take $\zeta^\param$ to be a standard family of mollifiers with respect to the 
Riemannian metric (for a precise definition, see \cite{LNO}). Here, we only record 
the properties which we will need, namely
\be
\label{10.47}
\int_M \xi^\param dV_{\gb} = 1,\quad |\nabla_\gb \xi^\param| \le \frac{C}{\param^{N+1}}, \quad 
|\Delta_\gb \xi^\param| \le \frac{C}{\param^{N+2}}.
\ee
In what follows we omit the superscript $\param$ when there is no risk of confusion.
One can check that the function $\ell_\gb$ satisfies
\[
d_q \ell_\gb (X)= \langle K_q, X\rangle_\gb  , \quad d_p \ell_\gb(X) = - \langle K_p, X\rangle_
\gb,
\]
for all vector fields $X$, where $K_p$ is the unit tangent vector at $p$ to the (unique) 
geodesic connecting $p$ and $q$. Therefore, for all vector fields $X$, we have
\be
\label{10.50}
d_q \xi(X) = \zeta'\circ \ell_\gb\langle K_q, X\rangle_\gb  , \quad d_p \xi(X) = - \zeta'\circ 
\ell_\gb\langle K_p, X\rangle_\gb.
\ee
Also, when there is no risk of confusion, we shall simply write $\zeta'$ instead of $
\zeta'\circ \ell_\gb$.
 We find from \eqref{10.45}
\[
\aligned 
&   \iint_{M \times M} \Big( - \xi d\chi\, F_p(u,v) - \chi \big(
d_p\xi\, F_p(u,v) +  d_q\xi\, F_q(u,v)\big)
\\
&\qquad \qquad -\chi  \xi \sgn(u-v)\big( \dive f_q(u) - \dive f_p(v)
\big)\Big) \,dV_g(p)\,dV_g(q) 
\\
& \le \iint_{M \times M} E(\chi \xi) \,dV_g(p)\,dV_g(q),
\endaligned
\]
where $u = u_p $ and $v=v_q$. We now write the last inequality as
\be
\label{10.60}
I_1 \le I_2 + I_3 + I_4,
\ee
with
\[
\aligned
&I_1 =\iint_{M \times M} - \xi d\chi\, F_p(u,v) \,dV_g(p)\,dV_g(q) ,
\\
&I_2 =\iint_{M \times M}  \chi \big( d_p\xi\, F_p(u,v) +  d_q\xi\, F_q(u,v)\big) \,dV_g(p)\,dV_g(q) ,
\\
&I_3 =\iint_{M \times M}  \chi  \xi \sgn(u-v)\big( \dive f_q(u) - \dive f_p(v) \big) \,dV_g(p)\,dV_g(q) ,
\\
&I_4 = \iint_{M \times M} E(\chi \xi )\,dV_g(p)\,dV_g(q).
\endaligned
\]

\subsection*{The term $I_2$}

Note that the term $I_2$ vanishes in the ``homogenous'' Euclidian case.
Writing $I_2 = \int_M \chi_p I'_2 dV_g(p)$, we can write
\be\label{10.65}
\aligned
I'_2 &= \int_M d_q \xi (F_q(u_p, v_p)) \,  dV_g(q)
 +  \int_M d_q \xi \Big( F_q(u_p, v_q) - F_q(u_p, v_p) \Big) dV_g(q)
\\
& \quad + \int_M d_p \xi \Big( F_p(u_p, v_q) - F_p(u_p, v_p) \Big) dV_g(q)
\\
&  = I_{2,1}+ I_{2,2} + I_{2,3},
\endaligned
\ee
since
\[
\int_M d_p \xi (F_p(u_p, v_p)) dV_g(q) = 0.
\]
To see this, note that $F_p(u_p, v_p)$ does not depend on $q$, and thus
\[
\aligned
F^i \int_M \del_{p^i} \xi dV(q) 
 = F^i \del_{p^i} \Big( \int_M  \xi dV(q) \Big) 
%\\
 = F^i \del_{p^i}(1) = 0.
\endaligned
\]
We have
\[
\aligned
\int_M d_q \xi (F_q(u_p, v_p)) dV_g(q) = - \int_M \xi \dive F_q(u_p, v_p) dV_g(q),
\endaligned
\]
since we can integrate by parts in this term, as $v$ and $u$ do not depend on the integration 
variable $q$. The remaining terms are estimated using the regularity of the flux and the 
difference $|v_p - v_q|$, as follows. From \eqref{10.65} we find
\be
\label{10.67}
\aligned
&I_{2,2} + I_{2,3} 
\\ 
&= \int_0^1 \int_M (v_q - v_p) \Big( d_q \xi \big( \del_v F_q(u_p , v^*) \big) 
 -  d_p \xi \big( \del_v F_p(u_p , v^*) \big) \Big) dV_g(q) ds,
\endaligned
\ee
with $v^* = sv_q + (1-s) v_p$.
From \eqref{10.05}, \eqref{10.50}, we have
\be
\label{10.70}
\aligned
& d_p \xi\, F_p(u,v) + d_q \xi\,  F_q(u,v)
= - \zeta' \big( \langle K_p , F_p(u,v)\rangle_\gb  + \langle K_q , F_q(u,v)   \rangle_\gb
\big). 
\endaligned
\ee

Now, 
consider the parallel transport of vectors along a curve on $M$ as
follows. Let $\gamma : [0,t_0] \to M$ be a smooth curve, and
$\nabla$ the covariant derivative operator associated to the
Riemannian metric $\gb$. Given $0\le s\le t\le t_0$, the parallel
transport is the operator $\tau_{s,t} : T_{\gamma(s)}M \to
T_{\gamma(t)} M$ such that, given a vector $ X_{\gamma(s)} \in
T_{\gamma(s)}M$, then $\tau_{s,t} X_{\gamma(s)} \in
T_{\gamma(t)}M$ is the unique solution of the differential
equation
\[
\nabla_{\gamma'} X =0, \qquad \tau_{s,t} X_{\gamma(s)}|_{t=s} = X_{\gamma(s)}.
\]
For our purposes, it is more important to note that, conversely, one may recover the 
covariant derivative from the notion of parallel transport, using the following relation,
\be
\label{10.80}
\nabla_{\gamma'(0)} X = \lim_{h\to 0} \frac{\tau_{h,0}X_{\gamma(h)} - X_{\gamma(0)}}
{h}.
\ee
Also, the parallel transport enjoys the property of preserving scalar products, that is, for all 
vector fields $X,Y$ defined along the curve $\gamma$,
\be
\label{10.90}
\langle X_{\gamma(s)},Y_{\gamma(s)} \rangle_{\gb_{\gamma(s)}} = \langle 
\tau_{s,t}X_{\gamma(s)}, \tau_{s,t}Y_{\gamma(s)} \rangle_{\gb_{\gamma(t)}}.
\ee
Furthermore, if $\gamma$ is a geodesic curve, then its tangent vector is invariant under 
parallel transport. With the notations above, if  $p=\gamma(0)$ and $q=\gamma(h)$, then
$ \tau_{0,h} K_p = K_q.$
Using the above properties of the parallel transport, we find 
\be
\label{10.95}
\aligned
\langle K_p , F_p(u,v)\rangle_\gb  - \langle K_q, F_q(u,v) \rangle_\gb 
&= -h \big\langle 
K_p , \frac{\tau_{h,0} F_q(u,v) - F_p(u,v)}{h} \big\rangle_\gb
\\
& =-h \langle K_{p^*} , \nabla_{K_{p^*}} F_{p^*}(u,v) \rangle_\gb 
\endaligned
\ee
with $h=\ell_\gb(p,q)$ and $p^*$ some point on the geodesic from $p$ to $q$.
We write simply $|X|$  for the norm of a vector with respect to the 
reference Riemannian metric to keep the notations simple. 
From \eqref{10.65},
\eqref{10.67} and \eqref{10.95}, with $\del_v F$ instead of $F$, we deduce 
\[
\aligned
I'_2 &\le \int_M \xi \dive F_q(u_p, v_p) dV_g(q)
\\
&\quad + \int_0^1 \int_{M}  h |v_p - v_q| \big|\big\langle \zeta' K_{p^*}, 
(\nabla_{K_{p^*}} \del_v F)_{p^*}(u_p,v^*) \big\rangle_\gb \big| 
\,dV_g(q)\, ds. 
\endaligned
\]
Now, since $|K| =1$ and $h\le\param$,
\[
\aligned
& \int_0^1 \int_{M}  h |v_p - v_q| \big|\big\langle \zeta' K_{p^*}, (\nabla_{K_{p^*}} \del_v F)_{p^*}
(u_p,v^*) \big\rangle_\gb \big| dV_g(q)\, ds
\\
&\le \Lambda_1\int_M 
\param |\zeta'| |v_p - v_q| dV_g(q),
\endaligned
\]
where $\Lambda_1$ is defined in \eqref{10.41}.
Thus
\be
\label{10.100}
\aligned
I_2 &\le - \iint_{M \times M}\chi \xi \dive F_q(u_p, v_p) dV_g(q)\,dV_g(p)
\\
& \quad +  \Lambda_1\iint_{M \times M} \chi \param |\zeta'| |v_p - v_q| dV_g(q)\,dV_g(p) .
\endaligned
\ee

\subsection*{The term $I_3$}

We now turn to the term $I_3$ in \eqref{10.60}, and write 
\[
\aligned
 \dive f_p(v_q) - \dive f_q(u_p)  & = \dive f_p (v_q) - \dive f_q (v_q)
 + \dive f_q (v_q) - \dive f_q (u_p), 
\endaligned
\]
thus
\[
\aligned
I_3 
= 
& -\iint_{M \times M}\chi \xi \sgn(u_p - v_q) \big\{ \dive f_p (v_q) 
- \dive f_q (v_q) \big\} dV_g(q)
\,dV_g(p)
\\
&- \iint_{M \times M}\chi \xi \sgn(u_p - v_q) \big\{ \dive f_q (v_q) 
 - \dive f_q (u_p) \big\} 
dV_g(q)\,dV_g(p).
\endaligned
\]
The first term is bounded by
\[
\iint_{M \times M}  \chi \xi  \param \Lambda_2 dV_g(q)\,dV_g(p),
\]
while the second term is simply
\[
\iint_{M \times M}\chi \xi \dive F_q (u_p, v_q) dV_g(q)\,dV_g(p).
\]
Thanks to \eqref{10.100}, this leads to
\[
\aligned
|I_2 + I_3| &\le \iint_{M \times M}  \chi \xi \bigl| \dive F_q(u_p,v_q) - \dive F_q(u_p, v_p)  
\bigr| dV_g(q)\,dV_g(p)
\\
& \quad + \iint_{M \times M} \Lambda_1 \chi \param |\zeta'| |v_p - v_q| + \param\Lambda_2\chi 
\xi \,dV_g(q)\,dV_g(p).
\endaligned
\]
We have
\[
\aligned
& \iint_{M \times M}  \chi \xi \bigl| \dive F_q(u_p,v_q) - \dive F_q (u_p, v_p)  \bigr| dV_g(q)\,dV_g(p)
\\
&  \le\iint_{M \times M} \chi \xi \Lip_u (\dive f) | v_p - v_q| dV_g(q)\,dV_g(p),
\endaligned
\]
which gives
\be
\label{10.110}
\aligned
|I_2 + I_3| \le \iint_{M \times M} \chi (\Lambda_1 \param |\zeta'|   + \Lambda_3 \xi )  |v_p - v_q | 
+ \param\Lambda_2\chi \xi \,dV_g(q)\,dV_g(p).
\endaligned
\ee

\subsection*{The term $I_1$}
We now treat the main term $I_1$. We take $\chi$ to be a function which is compactly 
supported in time and constant along the hypersurfaces $\Hcal_t$ and, thus, for all tangent 
vectors $Y$, we have $d\chi(Y) = \del_t\chi Y^t$. 
First, we have
\[
\aligned
I_1&= \iint_{M \times M} - \xi \del_t\chi  F_p^t(u_p,v_p)   \,dV_g(p)\,dV_g(q)
\\
& \quad - \iint_{M \times M} \xi \del_t\chi  \big( F_p^t(u_p,v_q) - F_p^t(u_p, v_p) \big)  \,dV_g(p)\,dV_g(q). 
\endaligned
\]
Now,
\[
\aligned
&\iint_{M \times M} \xi \del_t\chi \big( F_p^t(u_p,v_q) - F_p^t(u_p, v_p) \big)  \,dV_g(p)\,dV_g(q)
\\
&\le \iint_{M \times M} \Lip_u f^t \xi |\del_t\chi|  |v_p - v_q|   \,dV_g(p)\,dV_g(q),
\endaligned
\]
where we have used the fact that, since $\del_uf^t > 0,$ we have $| f^t(v) - f^t(u) | = 
F^t(u,v).$
From the last inequality and \eqref{10.100}, \eqref{10.110}, the inequality \eqref{10.60} 
becomes
\be
\label{10.120}
\aligned
& \iint_{M \times M} - \xi \del_t\chi \, F_p^t(u_p,v_p) \,dV_g(p)\,dV_g(q)
\\
&\le  \iint_{M \times M} \Phi_1(p,q)  |v_p - v_q|  \,dV_g(p)\,dV_g(q)
 +   \iint_{M \times M} \param \chi \xi \Lambda_2   \,dV_g(p)\,dV_g(q) + I_4,
\endaligned
\ee
with
\be
\label{10.130}
\aligned
&\Phi_1 =  \chi (\Lambda_1 \param |\zeta'| + \xi\Lambda_3) + |\del_t \chi| \xi
\Lambda_0
\endaligned
\ee
and where the constants $\Lambda_i$ are given by \eqref{10.41}.
Consider now the first term of \eqref{10.120}. From the properties of the mollifiers $\xi$, 
namely $\int_M \xi_{p,q} dV_g(q) = 1$ for all $p\in M$, we find
\be
\label{10.135}
\aligned
\iint_{M \times M} - \xi \del_t\chi \, F_p^t(u_p,v_p)  \,dV_g(p)\,dV_g(q) &= \int_M -\del_t\chi \, 
F_p^t(u_p,v_p)\,dV_g(p).
\endaligned
\ee
We choose the function $\chi = \chi^\eps \in(0,1)$ (which only depends on the $t$ 
coordinate and the small parameter $\eps> 0$), to be identically one if $t\in (\eps,T)$, 
supported on the set
\[
\bigcup_{t\in(0,T+\eps)}\Hcal_t ,
\]
so that its derivative, $\del_t\chi,$ is supported in
\[
\bigcup_{t\in(0,\eps)\cup(T,T+\eps)}\Hcal_t ,
\]
and satisfying
$
\chi^\eps \to \one_{(\Hcal_t )_{t\in(0,T)}}
$
as $\eps\to 0$. Therefore, since $\del_t \chi$ approaches $\delta_{\Hcal_0} - 
\delta_{\Hcal_T}$, and in view of the regularity assumptions on $u$ and $v$
and of 
\eqref{10.135}, we get 
\be
\label{10.140}
\aligned
& \limsup_{\eps\to 0} \int_M -\del_t\chi^\eps \, F_p^t(u_p,v_p)\,dV_g(p)
\\
&\le \int_{\Hcal_T} F_p^t(u_p,v_p) \,dV_g(p)
 - \int_{\Hcal_0} F_p^t(u_p,v_p) \,dV_g(p).
\endaligned
\ee

Next, consider the second term in \eqref{10.120}. It yields the term $E^\delta_v$ depending on the 
regularity of $v$, as follows. In view of \eqref{10.130} and since there exists a 
constant $C_N$ such that
\be
\label{10.150}
\aligned
\xi\Lambda_3 + \Lambda_1\param |\zeta'|  \le \frac{C_N}{\param^{N+1}}(\Lambda_1 + 
\Lambda_3), \qquad
\sup_{p\in M}\frac{|B_p(\param)|}{\param^{N+1}} \le C_N, 
\endaligned
\ee
we find
\[
\aligned
& \limsup_{\eps\to0}\Lambda_0 \iint_{M \times M} \xi  |\del_t \chi^\eps|  |v_p - v_q|  \,dV_g(q)\,dV_g(p) 
\\
&\le 
{C_N} \Lambda_0 \max_{t=0,T}  \int_{\Hcal_t} \dashint_{B_p(\param)} |v_p- v_q| \,dV_g(q)\,dV_g(p)
\endaligned
\]
and
\[
\aligned
&\limsup_{\eps\to0} \iint_{M \times M} \chi^\eps \, (\xi\Lambda_3 + \Lambda_1\param |\zeta'|)  |v_p - 
v_q|  \,dV_g(q)\,dV_g(p) 
\\
&\le TC_N(\Lambda_1 + \Lambda_3)\sup_{t\in0,T}  \int_{\Hcal_t} 
\dashint_{B_p(\param)} |v_p- v_q| \,dV_g(q)\,dV_g(p).
\endaligned
\]
Thus, we have 
\be
\label{10.160}
\aligned
&\limsup_{\eps\to0} \iint_{M \times M} \Phi_1(p,q)  |v_p - v_q|  \,dV_g(q)\,dV_g(p) 
\\
&\le C_N(T\Lambda_1 + T\Lambda_3 + \Lambda_0) \sup_{t\in0,T}  \int_{\Hcal_t} 
\dashint_{B_p(\param)} |v_p- v_q| \,dV_g(q)\,dV_g(p).
\endaligned
\ee

Now, considering the last term in \eqref{10.120}, we find 
\be
\iint_{M \times M} \param \Lambda_2\chi\xi \,dV_g(q)\,dV_g(p) =\param \Lambda_2\int_M \chi \,dV_g(p)
\le 
\label{10.170}
T\param\sup_{t\in (0,T)} |\Hcal_t| 
\Lambda_2.
\ee
Finally, consider the error term $I_4$ in \eqref{10.60}, given by \eqref{10.43}. First, note that
\[
d_p\varphi = \xi d \chi + \chi d_p\xi, \quad \nabla_{g^t} \varphi = \chi \nabla_{g^t} \xi, 
\quad \Delta_{g^t} \varphi  = \chi \Delta_{g^t} \xi. \]
From the properties of the test-functions $\xi$ and $\chi^\eps$, \eqref{10.47}, 
from \eqref{10.150}, and using the regularity of $\alpha_H$, we find
\[
\aligned
\limsup_{\eps\to 0} \iint_{M \times M} |d\chi^\eps | \xi \alpha_H\,dV_g(p)\,dV_g(q)
&\le \limsup_{\eps\to 0} 
\int_M |\del_t \chi^\eps| \alpha_H \,dV_g(p)
\\
&\le \sum_{t=0,T}\int_{\Hcal_t} \alpha_H \,dV_g(p),
\endaligned
\]
\[
\aligned
& \limsup_{\eps\to 0} \iint_{M \times M} 
|d_p\xi| \chi^\eps \, \alpha_H  \,dV_g(q)\,dV_g(p) 
\\
& \le \limsup_{\eps\to 0} \int_M \chi^\eps \int_{B_p(\param)} |\zeta'_{p,q} | \alpha_H \,dV_g(q)\,dV_g(p)
\\
&\le \limsup_{\eps\to 0} \int_M \chi^\eps\frac{1}{\param}\dashint_{B_p(\param)} 
\frac{| B_p(\param)|}{\delta^{N+1}} \,dV_g(q)\,dV_g(p)
\\
&\le \limsup_{\eps\to 0} 
\frac{C}{\param}\int_M \chi^\eps \alpha_H \,dV_g(p)
 \le \frac{C}{\param} \int_{M_T} \alpha_H \,dV_g(p),
\endaligned
\]
then 
\[
\aligned
\limsup_{\eps\to 0} \iint_{M \times M} \xi \chi^\eps \alpha_K  \,dV_g(p)\,dV_g(q) & \le\int_{M_T} \alpha_K \,dV_g(p),
\endaligned
\]
and
\[
\aligned
& \limsup_{\eps\to 0}\iint_{M \times M}  \alpha_L \alpha_a (  |\Delta_{g^t} \varphi| + |\nabla_{g^t} \varphi|) 
\,dV_g(p)\,dV_g(q)
\\
&\le C\big( \frac{1}{\param} + \frac{1}{\param^2} \big)\int_{M_T} \alpha_L\alpha_a \,dV_g(p).
\endaligned
\]
The estimate \eqref{10.42} now follows from the 
inequalities above and from \eqref{10.140}, 
\eqref{10.160}, \eqref{10.170}. This completes the proof of Theorem~\ref{10-10}.

%===========================================================================================================

\section*{Acknowledgements}

This paper was written when the second author (PLF) participated in the international research program on 
``Nonlinear Partial Differential Equations'',
 hold at the Centre for Advanced Study,
  Norwegian Academy of Sciences and Letters during the Academic Year 2008--09. This author is grateful to Helge Holden and 
Kenneth Karlsen for their invitation and hospitality. 
 
The authors were supported by the Agence Nationale de la Recherche (ANR) via the grant 06-2-134423. 
PA was also supported by the Portuguese Foundation for Science and Technology (FCT) 
through the post-doctoral fellowship SFRH/\-BPD/\-43548/\-2008, and by a Ci\^encia 2008 fellowship.
WN was also supported by FAPERJ via the grant
E-26/ 111.564/2008 entitled ``Analysis, Geometry and Applications'',
 by Pronex-FAPERJ through the grant E-26/ 110.560/2010 entitled ``Nonlinear Partial Differential Equations'', and by a joint project between Brazil and France.

%======================================================================================================================
\small 
\newcommand{\auth}{\textsc}

\end{document}